\documentclass[11pt]{amsart}

\usepackage{amsfonts}
\usepackage{amssymb}
\makeatletter \@addtoreset{equation}{section}

\makeatother \makeatletter

\newcommand{\bea}{\begin{eqnarray}}
\newcommand{\eea}{\end{eqnarray}}
\def \1{\bf 1}
\def \Zset{\Bbb Z}
\def \Cset{\Bbb C}

\def \Nset{\Bbb N}

\newtheorem{theorem}{Theorem}[section]
\newtheorem{corollary}[theorem]{Corollary}
\newtheorem{lemma}[theorem]{Lemma}
\newtheorem{proposition}[theorem]{Proposition}
\theoremstyle{definition}

\theoremstyle{remark}
\newtheorem{remark}[theorem]{Remark}

\begin{document}

\title{ Representations of certain non-rational vertex
operator algebras of affine type}

\author{Dra\v{z}en Adamovi\' c and  Ozren Per\v{s}e}
\address{Department of Mathematics, University of Zagreb, Croatia}
\email{adamovic@math.hr} \email{perse@math.hr}

\maketitle
\begin{abstract}
In this paper we study a series of vertex operator algebras of
integer level associated to the affine Lie algebra $A_{\ell}^{(1)}$.
These vertex operator algebras are constructed by using the explicit
construction of certain singular vectors in the universal affine
vertex operator algebra $N_{l}(n-2,0)$ at the integer level. In the
case $n=1$ or $l=2$, we explicitly determine Zhu's algebras and
classify all irreducible modules in the category $\mathcal{O}$. In
the case $l=2$, we show that the vertex operator algebra
$N_2(n-2,0)$ contains two linearly independent singular vectors of
the same conformal weight.
\end{abstract}

\section{Introduction}

Let ${\hat {\frak g}}$ be the affine Lie algebra associated with the
finite dimensional simple Lie algebra ${\frak g}$.  Then, on the
generalized Verma module $N(k,0)$,  exists the natural structure of
a vertex operator algebra  (cf. \cite{FrB,FZ,K2,LL,L,MP}). In the
representation theory of affine Lie algebras it is important
 to study annihilating ideals of highest weight representations.
 This problem is related with the ideal lattice of the vertex
 operator algebra $N(k,0)$ (cf. \cite{FM}).
 Every ${\hat {\frak g}}$--submodule $I$ of $N(k,0)$ becomes an ideal in
the vertex operator algebra $N(k,0)$, and on the quotient $N(k,0) /
I $, there exists the structure of a vertex operator algebra. Let
$N^{1}(k,0)$ be the maximal ideal in $N(k,0)$. Then  $L(k,0) =
N(k,0) / N^{1}(k,0)$ is a simple vertex operator algebra.

It is particularly important to study ideals in $N(k,0)$ generated
by one singular vector. If $k$ is a positive integer, then
$N^{1}(k,0)$ is generated by the singular vector $e_\theta(-1)
^{k+1} {\bf 1}$ (cf. \cite{FZ,MP}). In the case when $k$ is an
admissible rational number, the maximal ideal is also generated by
one singular vector, but the expression for this singular vector is
more complicated (cf. \cite{A1,AM,P}).

By using explicit formulas for singular vectors which involve the
powers of determinants, a family of ideals were constructed in
\cite{A4}. The corresponding quotient vertex operator algebras give
new examples of vertex operator algebras of affine type which are
different from $N(k,0)$ and $L(k,0)$. So it is an interesting
problem to investigate their representation theory.

In this paper we shall study the (non-simple) vertex operator
algebras associated to the affine Lie algebra $A_{l} ^{(1)}$ of
integer level. We shall classify all irreducible highest weight
representations of the  vertex operator algebra $V_{l,1}$ associated
to a level $-1$ representation for $A_l ^{(1)}$. We use the methods
developed in \cite{A1} and \cite{P}. We also classify the
finite-dimensional irreducible representations of corresponding
Zhu's algebra $A(V_{l,1})$. We prove  an interesting fact that every
finite-dimensional irreducible $A(V_{l,1})$--module has all
$1$--dimensional weight spaces. Therefore Zhu's algebra $A(V_{l,1})$
provides a natural framework for studying irreducible
representations having $1$--dimensional weight spaces. These
representations appeared in \cite{BL} in the context of the
representation theory of Lie algebras. So our results show their
importance in the theory of vertex operator algebras.

We shall also extend the construction of singular vectors   from
\cite{A4}.  In particular, we shall present new explicit formulas
for singular vectors in the case of affine Lie algebra $A_2 ^{(1)}$
at integer levels. We show that the vertex operator algebra
$N_2(n-2,0)$ contains two linearly independent  singular vectors
$v_{2,n}$ and $\widetilde{v_{2,n}}$ of conformal weight $3n$.
Moreover, there is an automorphism $\Psi$ of the vertex operator
algebra $N_2(n-2,0)$ such that $\Psi(v_{2,n}) =
\widetilde{v_{2,n}}$.
 These singular vectors generate
ideals $J_2(n-2,0)$ and $\widetilde{J_2} (n-2,0)$ whose top levels
are irreducible modules with all $1$--dimensional weight spaces. Let
$V_{2,n}$ and $\widetilde{V_{2,n}}$ be the corresponding quotient
vertex operator algebras. It turns out that these vertex operator
algebras are isomorphic. We show  the irreducible representations of
$V_{2,n}$  and $\widetilde{V_{2,n}}$
 in the category $\mathcal{O}$ are parameterized by the zeros of the
same polynomial $P \in {\Cset}[x,y]$.

Therefore, for every $n \in {\Nset}$, we have   a new  vertex
operator algebra, which is a certain quotient of $N_{2}(n-2,0)$, for
which we know all irreducible highest weight modules. The
irreducible highest weight modules for that vertex operator algebra
describe  a new infinite class of $A_2 ^{(1)}$--modules, which
includes integrable highest weight $A_2 ^{(1)}$--modules of level
$n-2$.

We believe that the categories of representations which appear in
this paper have interesting characters, and we hope to study the
corresponding fusion rules.

In this paper, ${\Zset}_{+}$ and $\Nset$ denote the sets of
nonnegative integers and positive integers, respectively.

\section{Vertex operator algebras associated to affine Lie algebras}

In this section we review certain results about vertex operator
algebras associated to affine Lie algebras.

Let $(V, Y, {\bf 1}, \omega)$ be a vertex operator algebra (cf.
\cite{B,FHL,FLM}). We shall assume that
$$V= \bigsqcup_{n \in {\Zset}_{+} } V_n, \quad \mbox{where} \ V_n =
\{  a \in V  \vert L(0) a = n a \}. $$
For $a \in V_n$ we shall write $ {\rm wt} (a)  = n$.
 Following \cite{Z}, we
define bilinear maps $* : V \times V \to V$ and $\circ  :  V \times
V \to V$ as follows. For any homogeneous $a \in V$ and for any $b
\in V$, let
\begin{eqnarray*}
a \circ b={\rm Res}_{z}\frac{(1+z)^{{\rm wt} (a)}}{z^{2}}Y(a,z)b  \\
a * b={\rm Res}_{z}\frac{(1+z)^{{\rm wt} (a)}}{z}Y(a,z)b
\end{eqnarray*}
and extend to $V \times V \to V$ by linearity. Denote by $O(V)$ the
linear span of elements of the form $a \circ b$, and by $A(V)$ the
quotient space $V/O(V)$. The multiplication $*$ induces the
multiplication on $A(V)$ and $A(V)$ has a structure of an
associative algebra. There is one-to-one correspondence between the
equivalence classes of the irreducible $A(V)$--modules and the
equivalence classes of the irreducible ${\Zset}_{+}$--graded weak
$V$--modules (\cite{Z}).

Let ${\frak g}$ be a simple Lie algebra over ${\Cset}$ with a
triangular decomposition  ${\frak g}={\frak n}_{-} \oplus {\frak h}
\oplus {\frak n}_{+}$. The affine Lie algebra $\hat{\frak g}$
associated to ${\frak g}$ is the vector space ${\frak g}\otimes
{\Cset}[t, t^{-1}] \oplus {\Cset}c $ equipped with the usual bracket
operation and the canonical central element~$c$ (cf. \cite{K}). Let
$h^{\vee}$ be the dual Coxeter number of $\hat{\frak g}$. Let
$\hat{\frak g}=\hat{\frak n}_{-} \oplus \hat{\frak h} \oplus
\hat{\frak n}_{+}$ be the corresponding triangular decomposition of
$\hat{\frak g}$.

Let $U$ be a ${\frak g}$--module, and let $k\in {\Cset}$. Let
$\hat{\frak g}_{+}={\frak g}\otimes t{\Cset}[t]$ act trivially on
$U$ and $c$ as the scalar multiplication operator $k$. Considering
$U$ as a ${\frak g}\oplus {\Cset}c \oplus \hat{\frak
g}_{+}$--module, we have the induced $\hat{\frak g}$--module (so
called {\it generalized Verma module})
\[
N(k,U)=U(\hat{\frak g})\otimes_{U({\frak g}\oplus {\Cset}c \oplus
\hat{\frak g}_{+})} U.
\]

For a fixed $\mu \in {\frak h}^{*}$, denote by $V(\mu)$ the
irreducible highest-weight ${\frak g}$--module with highest weight
$\mu$. We shall use the notation $N(k, \mu)$ to denote the
$\hat{\frak g}$--module $N(k,V(\mu))$. Denote by $N^1(k, \mu)$ the
maximal proper submodule of $N(k, \mu)$ and by $L(k, \mu)=N(k, \mu)
/N^1(k, \mu)$ the corresponding irreducible $\hat{\frak g}$--module.

For $k \neq - h^{\vee}$, $N(k,0)$ has the structure of vertex
operator algebra (cf. \cite{FZ} and \cite{MP}, see also
\cite{FrB,K2,LL,L}). The associative algebra $A(N(k,0))$ is
canonically isomorphic to $U (\frak g ) $ (\cite{FZ}).

Let $J=U(\hat{\frak g})v$ be the  $\hat{\frak g}$--submodule of
$N(k,0)$ generated by the  singular vector $v$. Then $J$ is an ideal
in $N(k,0)$ and the quotient $N(k,0) /J$ is also a vertex operator
algebra.  We present the method from \cite{A1,AM,MP,P} for the
classification of irreducible $A(N(k,0)/J)$--modules from the
category $\mathcal{O}$ by solving certain systems of polynomial
equations.

Let $v'$ be the image of the vector $v$ in Zhu's algebra $A(N(k,0))
\cong U (\frak g )$. Denote by $_L$ the adjoint action of  $U(\frak
g)$ on $U(\frak g)$ defined by $ X_Lf=[X,f]$ for $X \in \frak g$ and
$f \in U(\frak g)$. Let $R$ be a $U(\frak g)$--submodule of $U(\frak
g)$ generated by the vector $v'$ under the adjoint action. Clearly,
$R$ is an irreducible highest-weight $U(\frak g)$--module. Let
$R_{0}$ be the zero-weight subspace of $R$.

The next proposition follows from \cite{A1,AM,MP,P}:

\begin{proposition} \label{prop-R0}
Let $V(\mu)$ be an irreducible highest weight $U(\frak g)$--module
with the highest weight vector $v_{\mu}$, for $\mu \in {\frak
h}^{*}$. The following statements are equivalent:
\item[(1)] $V(\mu)$ is an $A(N(k,0)/J)$--module,
\item[(2)] $RV(\mu)=0$,
\item[(3)] $R_{0}v_{\mu}=0.$
\end{proposition}

Let $r \in R_{0}$. Clearly there exists the unique polynomial $p_{r}
\in S( \frak h)$ such that
\[
rv_{\mu}=p_{r}(\mu)v_{\mu}.
\]
Set $ {\mathcal P}_{0}=\{ \ p_{r} \ \vert \ r \in R_{0} \}.$ We
have:

\begin{corollary} \label{c.0.8} There is one-to-one correspondence between
\item[(1)] irreducible $A(N(k,0)/J)$--modules
from the category $\mathcal{O}$,
\item[(2)] weights $\mu \in {\frak h}^{*}$ such that
$p(\mu)=0$ for all $p \in {\mathcal P}_{0}$.
\end{corollary}

\section{Simple Lie algebra of type $A_{l}$}

Let $\frak g$ be the simple Lie algebra of type $A_{l}$. The root
system of $\frak g$ is given by
\[
\Delta=\{\epsilon_i - \epsilon_j \ \vert \ 1 \leq i,j \leq l+1, i
\neq j \}
\]
with $\alpha_1=\epsilon_1-\epsilon_2,...,\alpha_{l}=\epsilon_{l}-
\epsilon_{l+1}$ being a set of simple roots. The highest root is
$\theta=\epsilon_{1}- \epsilon_{l+1}$. We fix the root vectors and
coroots as in \cite{A4} and \cite{FF}.  For any positive root
$\alpha \in \Delta_{+}$ denote by $e_{\alpha}$ and $f_{\alpha}$ the
root vectors corresponding to $\alpha$ and $- \alpha$, respectively.
Denote by $h_{\alpha}$ the corresponding coroots, and by
$h_{i}=h_{\alpha_{i}}$, for $i =1, \dots, l$. Let $\omega_1, \dots ,
\omega_l$ be the fundamental weights for $A_{l}$.

We recall the following well-known and important fact about the
weight spaces of $\frak g$-modules $V(n \omega_1)$ and $V(n
\omega_l)$.

\begin{proposition} \label{one-dim}
All weight spaces of irreducible $\frak g$--modules $V(n \omega_1)$
and $V(n \omega_l)$ are $1$--dimensional.
\end{proposition}

\begin{remark} One can show that every irreducible finite-dimensional
$\frak g$--module having all $1$--dimensional weight spaces is
isomorphic to either $V(n \omega_1)$ or $V(n \omega_l)$ for certain
$n \in {\Nset}$. In what follows we shall see that the category of
irreducible modules with $1$--dimensional weight spaces provides all
irreducible finite-dimensional modules for certain Zhu's algebra.
\end{remark}

\section{Singular vectors in the vertex operator algebras
$N_{l}(n-2,0)$ for $n \in {\Nset}$} \label{sec-a3}

Let $\hat{\frak g}$ be the affine Lie algebra of type $A_{l}^{(1)}$.
For $k\in \Cset$, denote by $N_{l}(k,0)$ the generalized Verma
module associated to $\hat{\frak g}$ of level $k$.

We shall recall some facts about singular vectors in the vertex
operator algebra $N_{l}(n-2,0)$ for $n \in {\Nset}$. We know that
$N_{l}(n-2,0)$ is not a simple vertex operator algebra. If $n \ge
2$,  it contains the maximal ideal $N_{l} ^{1} (n-2,0)$ which is
generated by the singular vector $e_{\theta} (-1) ^{n-1} {\bf 1}$.
In the case $l=1$, one can show that the maximal ideal is an
irreducible $\hat{\frak sl_2}$--module. Therefore, there are no
other ideals in $N_{l} ^{1} (n-2,0)$. But for $l\ge 2$, the
situation is very much different. We have non-trivial ideals which
are different from $N_{l} ^{1} (n-2,0)$. These ideals are generated
by singular vectors in $N_{l} ^{1} (n-2,0)$. Some non-trivial ideals
were constructed in \cite{A4} in the case $l \ge 3$.

We have:

\begin{theorem} [{\cite[Theorem 5.1]{A4}}] \label{sing-a3}
Assume that $l, n \in {\Nset}$, $l \ge 3$.
\item[(1)]For every $n \in {\Nset}$, vector
\begin{eqnarray*}
v_{l,n}=\left( e_{\epsilon_1 - \epsilon_{l+1}}(-1)e_{\epsilon_2 -
\epsilon_{l}}(-1)  - e_{\epsilon_2 -
\epsilon_{l+1}}(-1)e_{\epsilon_1 - \epsilon_{l}}(-1)\right) ^{n
}{\bf 1}
\end{eqnarray*}
is a non-trivial singular vector in $N_{l}(n-2,0)$.
\item[(2)] Assume that $n \ge 2$. Then the maximal submodule
$N_{l} ^{1} (n-2,0)$ of $N_{l}(n-2,0)$ is reducible.
\end{theorem}

\begin{remark}
Note that we constructed singular vectors in the case $l \ge 3$. So
the construction from \cite{A4} can not be applied directly in the
case $l=2$. In \cite{FF}, A. Feingold and I. Frenkel presented an
explicit realization of the vacuum $A_{l}^{(1)}$--module of level
$-1$. Let us denote this module by $\widetilde{L}_{l} (-1,0)$.
$\widetilde{L}_{l} (-1,0)$ carries the structure of a vertex
operator algebra (this vertex operator algebras was also studied in
\cite{KR}), and it is a non-trivial quotient of the vertex operator
algebra $N_{l}(-1,0)$. By  a simple analysis of characters, one can
show that in the case $l \ge 3$, $v_{\ell,1}$ is a unique, up to
scalar factor, singular vector in $N_{l}(-1,0)$ with conformal
weight $2$. We omit details. In the case $l=2$, we have:
\begin{eqnarray}
\mbox{ch}\widetilde{L}_{2} (-1,0) = 1 + 8 q + 44 q^{2} + 172 q^{3}+
o(q^{4});  \nonumber \\ \mbox{ch}N_{2} (-1,0) = 1 + 8 q + 44 q^{2} +
192 q^{3}+ o(q^{4}) \label{char-nej}
\end{eqnarray}
 which  indicates that there are no singular
vectors of conformal weight $2$, and that there should exist at
least one singular vector of conformal weight $3$ in  $N_2(-1,0)$.
In Section \ref{sl-3} we shall find two linearly independent
singular vectors of conformal weight $3$. Moreover, we shall extend
Theorem \ref{sing-a3} in the case $l=2$.
\end{remark}

\section{Vertex operator algebras associated to affine Lie algebras of
type $A_{l}^{(1)}$ with level $-1$, for $l \geq 3$}

In the special case $n=1$, Theorem \ref{sing-a3} implies that vector
\begin{eqnarray}
&& v_{l,1}=e_{\epsilon_1 - \epsilon_{l+1}}(-1)e_{\epsilon_2 -
\epsilon_{l}}(-1) {\bf 1} - e_{\epsilon_2 -
\epsilon_{l+1}}(-1)e_{\epsilon_1 - \epsilon_{l}}(-1) {\bf 1}
\label{singul-l}
\end{eqnarray}
is a singular vector in $N_{l}(-1,0)$.

Define the ideal $J_{l,1}(-1,0)$ in the vertex operator algebra
$N_{l}(-1,0)$ by $ J_{l,1}(-1,0)=U(\hat{\frak g})v_{l,1}, $ and
denote the corresponding quotient vertex operator algebra by
$V_{l,1}$, i.e., $V_{l,1}=N_{l}(-1,0)/J_{l,1}(-1,0)$.

 Then
\cite[Theorem 5.4]{A4}:
\begin{proposition}
The associative algebra $A(V_{l,1})$ is isomorphic to the algebra
$U(\frak g)/I_{l,1}$, where $I_{l,1}$ is the two-sided ideal of
$U(\frak g)$ generated by
\begin{eqnarray*}
v' _{l,1} =e_{\epsilon_1 - \epsilon_{l+1}}e_{\epsilon_2 -
\epsilon_{l}} - e_{\epsilon_2 - \epsilon_{l+1}}e_{\epsilon_1 -
\epsilon_{l}}.
\end{eqnarray*}
\end{proposition}

The $U(\frak g)$--submodule $R$ of $U(\frak g)$ generated by the
vector $v'_{l,1}$ under the adjoint action is isomorphic to
$V(\omega_{2}+\omega_{l-1})$. In the next lemma we determine the
dimension of $R_{0}$.

\begin{lemma} \label{l.1.8}
\[
\dim R_{0} = \frac{(l-2)(l+1)}{2}
\]
\end{lemma}
{\bf Proof:} In this proof we use induction on $l$. We use the
notation $V_{l}(\mu)$ for the highest weight module for simple Lie
algebra of type $A_{l}$, with the highest weight $\mu \in {\frak
h}^{*}$.

For $l=3$ it can be easily checked that $ \dim R =20$ and $ \dim
R_{0} =2$. Suppose that the claim of this lemma holds for simple Lie
algebra of type $A_{l-1}$, $l-1 \geq 3$. Let $\frak g$ be simple Lie
algebra of type $A_{l}$. Let $\frak g'$ be the subalgebra of $\frak
g$ associated to roots $\alpha_{1}, \ldots ,\alpha_{l-1}$. Then
$\frak g'$ is a simple Lie algebra of type $A_{l-1}$. We can
decompose $\frak g$--module $V_{l}(\omega_{2}+\omega_{l-1})$ into
the direct sum of irreducible $\frak g'$--modules. If we denote by
$v$ the highest weight vector of $\frak g$--module
$V_{l}(\omega_{2}+\omega_{l-1})$, then it can be easily checked that
$v$, $f_{\epsilon_{l-1} - \epsilon_{l+1}}.v$, $(2 f_{\epsilon_{2} -
\epsilon_{l+1}}- f_{\epsilon_{2} - \epsilon_{3}}f_{\epsilon_{3} -
\epsilon_{l+1}}).v$ and $f_{\epsilon_{l-1} - \epsilon_{l+1}}(2
f_{\epsilon_{2} - \epsilon_{l+1}}- f_{\epsilon_{2} -
\epsilon_{3}}f_{\epsilon_{3} - \epsilon_{l+1}}).v$ are highest
weight vectors for $\frak g'$, which generate $\frak g'$--modules
isomorphic to $V_{l-1}(\omega_{2}+\omega_{l-1})$,
$V_{l-1}(\omega_{2}+\omega_{l-2})$,
$V_{l-1}(\omega_{1}+\omega_{l-1})$ and
$V_{l-1}(\omega_{1}+\omega_{l-2})$, respectively. It follows from
Weyl's formula for the dimension of irreducible module that $ \dim
V_{l}(\omega_{2}+\omega_{l-1})=\frac{(l-2)(l+2)(l+1)^2}{4}$, $ \dim
V_{l-1}(\omega_{2}+\omega_{l-1})= \dim
V_{l-1}(\omega_{1}+\omega_{l-2})= \frac{(l-2)l(l+1)}{2}$, and $ \dim
V_{l-1}(\omega_{1}+\omega_{l-1}) =(l-1)(l+1)$, since
$V_{l-1}(\omega_{1}+\omega_{l-1})$ is the adjoint module for $\frak
g'$. We obtain
\begin{eqnarray*}
&& V_{l}(\omega_{2}+\omega_{l-1}) \cong
V_{l-1}(\omega_{2}+\omega_{l-1}) \oplus
V_{l-1}(\omega_{2}+\omega_{l-2}) \oplus
V_{l-1}(\omega_{1}+\omega_{l-1}) \\
&& \quad \oplus V_{l-1}(\omega_{1}+\omega_{l-2}).
\end{eqnarray*}
Clearly, there are no zero-weight vectors for $\frak g$ in $\frak
g'$--modules generated by highest weight vectors $v$ and
$f_{\epsilon_{l-1} - \epsilon_{l+1}}(2 f_{\epsilon_{2} -
\epsilon_{l+1}}- f_{\epsilon_{2} - \epsilon_{3}}f_{\epsilon_{3} -
\epsilon_{l+1}}).v$. There are $l-1$ linearly independent
zero-weight vectors for $\frak g'$ in the module
$V_{l-1}(\omega_{1}+\omega_{l-1})$, since it is isomorphic to the
adjoint module for $\frak g'$, and all those vectors have weight
zero for $\frak g$. The inductive assumption implies that there are
$\frac{(l-3)l}{2}$ linearly independent zero-weight vectors for
$\frak g$ in $\frak g'$--module $V_{l-1}(\omega_{2}+\omega_{l-2})$,
which implies $ \dim R_{0} = \frac{(l-3)l}{2}+(l-1)=
\frac{(l-2)(l+1)}{2}$. \qed

\begin{lemma} \label{l.1.9}
Let
\begin{eqnarray*}
& & (1) \ p_{ij}(h)=h_{i}h_{j}, \
\mbox{for } i=1, \ldots ,l-2, j-i \geq 2 \\
& & (2) \ q_{i}(h)=h_{i}(h_{i-1}+h_{i}+h_{i+1}+1),  \ \mbox{for }
i=2, \ldots ,l-1.
\end{eqnarray*}
Then $p_{ij},q_{i} \in {\mathcal P}_{0}$.
\end{lemma}
{\bf Proof:} We use induction on $l$. For $l=3$, one can easily
obtain that
\begin{eqnarray*}
&& (-f_{\epsilon_1 - \epsilon_4}f_{\epsilon_2 - \epsilon_3}-
f_{\epsilon_1 - \epsilon_3}f_{\epsilon_2 - \epsilon_4})_L v_{3,1} '
\in h_{1}h_{3} + U(\frak g){\frak n}_{+} \quad \mbox{and} \\
&& (f_{\epsilon_1 - \epsilon_4}f_{\epsilon_2 - \epsilon_3})_L
v_{3,1} ' \in h_{2}(h_{1}+h_{2}+h_{3}+1) + U(\frak g){\frak n}_{+},
\end{eqnarray*}
which implies that $p_{13},q_{2} \in {\mathcal P}_{0}$. Suppose that
the claim of this lemma holds for simple Lie algebra of type
$A_{l-1}$, $l-1 \geq 3$. Let $\frak g$ be simple Lie algebra of type
$A_{l}$. Let $\frak g'$ be the subalgebra of $\frak g$ associated to
roots $\alpha_{1}, \ldots ,\alpha_{l-1}$, and $\frak g''$ the
subalgebra of $\frak g$ associated to roots $\alpha_{2}, \ldots
,\alpha_{l}$. Then $\frak g'$ and $\frak g''$ are simple Lie
algebras of type $A_{l-1}$. Since
\begin{eqnarray*}
&& (f_{\epsilon_1 - \epsilon_3})_L v_{l,1} ' =e_{\epsilon_2 -
\epsilon_{l+1}}e_{\epsilon_3 - \epsilon_{l}} - e_{\epsilon_3 -
\epsilon_{l+1}}e_{\epsilon_2 -
\epsilon_{l}} \quad \mbox{and} \\
&& (-f_{\epsilon_{l-1} - \epsilon_{l+1}})_L v_{l,1} ' =e_{\epsilon_1
- \epsilon_{l}}e_{\epsilon_2 - \epsilon_{l-1}} - e_{\epsilon_1 -
\epsilon_{l-1}}e_{\epsilon_2 - \epsilon_{l}}=v_{l-1, 1}',
\end{eqnarray*}
we get the corresponding vectors for $\frak g''$ and $\frak g'$,
respectively. Using the inductive assumption for $\frak g'$, we get
that $p_{ij} \in {\mathcal P}_{0}$ for $i=1, \ldots ,l-3, j-i \geq
2$, and $q_{i} \in {\mathcal P}_{0}$ for $i=2, \ldots ,l-2$. And
using the inductive assumption for $\frak g''$, we get that $p_{ij}
\in {\mathcal P}_{0}$ for $i=2, \ldots ,l-2, j-i \geq 2$, and $q_{i}
\in {\mathcal P}_{0}$ for $i=3, \ldots ,l-1$. Also, one can easily
verify that
\begin{eqnarray*}
(-f_{\epsilon_1 - \epsilon_{l+1}}f_{\epsilon_2 - \epsilon_{l}}-
f_{\epsilon_2 - \epsilon_{l+1}}f_{\epsilon_1 - \epsilon_{l}})_L
v_{l,1} ' \in h_{1}h_{l} + U(\frak g){\frak n}_{+},
\end{eqnarray*}
which implies that $p_{1l} \in {\mathcal P}_{0}$. Thus $p_{ij} \in
{\mathcal P}_{0}$ for $i=1, \ldots ,l-2, j-i \geq 2$, and $q_{i} \in
{\mathcal P}_{0}$ for $i=2, \ldots ,l-1$, and the claim of lemma is
proved. \qed

\begin{lemma}
Polynomials $p_{ij} \in {\mathcal P}_{0}$ for $i=1, \ldots ,l-2, j-i
\geq 2$ and $q_{i} \in {\mathcal P}_{0}$ for $i=2, \ldots ,l-1$ form
a basis for the vector space ${\mathcal P}_{0}$.
\end{lemma}
{\bf Proof:} Lemma \ref{l.1.9} implies that $p_{ij} \in {\mathcal
P}_{0}$ for $i=1, \ldots ,l-2, j-i \geq 2$ and $q_{i} \in {\mathcal
P}_{0}$ for $i=2, \ldots ,l-1$ are $\frac{(l-2)(l+1)}{2}$ linearly
independent polynomials in the set ${\mathcal P}_{0}$. It follows
from Lemma \ref{l.1.8} that $\dim {\mathcal P}_{0} =
\frac{(l-2)(l+1)}{2}$. Thus, these polynomials form a basis for
${\mathcal P}_{0}$. \qed

\begin{proposition} The set
\begin{eqnarray*}
\left\{ V(t \omega _{1}) \ \vert \ t \in \Cset \right\} \cup \left\{
V(t \omega _{l}) \ \vert \ t \in \Cset \right\} \cup
\bigcup_{i=1}^{l-1} \left\{ V(t \omega _{i} + (-1-t) \omega _{i+1})
\ \vert \ t \in \Cset \right\}
\end{eqnarray*}
provides the complete list of irreducible $A(V_{l,1})$--modules from
the category $\mathcal{O}$.
\end{proposition}
{\bf Proof:} It follows from Corollary \ref{c.0.8} that the highest
weights $\mu \in {\frak h}^{*}$ of irreducible $A(V_{l,1})$--modules
$V(\mu)$ are in one-to-one correspondence with solutions of the
system of polynomial equations
\begin{eqnarray}
& &  \ h_{i}h_{j}=0, \ \mbox{for } i=1, \ldots ,l-2, j-i \geq 2; \label{1.1} \\
& &  \ h_{i}(h_{i-1}+h_{i}+h_{i+1}+1)=0,  \ \mbox{for } i=2, \ldots
,l-1. \label{1.2}
\end{eqnarray}
Clearly, $h_{i}=0$, for $i=1, \ldots ,l$ is a solution of the system
corresponding to the highest weight $\mu =0$. Suppose that $h_{i}
\neq 0$, for some $i \in \{2, \ldots ,l-1 \}$. Then (\ref{1.1})
implies that $h_{j}=0$, for $j \neq i-1,i,i+1$. From relations
\begin{eqnarray*}
& &  \ h_{i-1}h_{i+1}=0,   \\
& &  \ h_{i-1}+h_{i}+h_{i+1}+1=0,
\end{eqnarray*}
we get the solutions corresponding to highest weights $\mu =t \omega
_{i-1} + (-1-t) \omega _{i}$ and $\mu =t \omega _{i} + (-1-t) \omega
_{i+1}$. If $h_{1} \neq 0$, then (\ref{1.1}) implies that $h_{j}=0$,
for $j =3, \ldots l$. Relation (\ref{1.2}) then implies that
\begin{eqnarray*}
& &  \ h_{2}(h_{1}+h_{2}+1)=0,
\end{eqnarray*}
and we get the solutions corresponding to highest weights $\mu =t
\omega _{1}$ and $\mu =t \omega _{1} + (-1-t) \omega _{2}$.
Similarly, if $h_{l} \neq 0$, we get the solutions corresponding to
highest weights $\mu =t \omega _{l}$ and $\mu =t \omega _{l-1} +
(-1-t) \omega _{l}$. \qed

\begin{corollary}
 The set
\begin{eqnarray*}
\left\{ V(t \omega _{1}) \ \vert \ t \in {\Zset}_{+} \right\} \cup
\left\{ V(t \omega _{l}) \ \vert \ t \in {\Zset}_{+} \right\}
\end{eqnarray*}
provides the complete list of irreducible finite-dimensional
$A(V_{l,1})$--modules. Moreover, every irreducible
finite-dimensional $A(V_{l,1})$--module has all  $1$--dimensional
weight spaces for $\frak h$.
\end{corollary}

It follows from Zhu's theory that:

\begin{theorem} \label{t.klas.nivo-1}
The set
\begin{eqnarray*}
\left\{ L(-1, t \omega _{1}) \ \vert \ t \in \Cset  \right\} \cup
\left\{ L(-1, t \omega _{l}) \ \vert \ t \in \Cset  \right\} \cup
\bigcup_{i=1}^{l-1} \left\{ L(-1, t \omega _{i}+ (-1-t) \omega
_{i+1}) \ \vert \ t \in \Cset  \right\}
\end{eqnarray*}
provides the complete list of irreducible weak $V_{l,1}$--modules
from the category $\mathcal{O}$.
\end{theorem}

\begin{remark} \label{rem.not-semi}
Theorem \ref{t.klas.nivo-1} gives the classification of irreducible
objects in the category of weak $V_{l,1}$--modules that are in the
category $\mathcal{O}$ as $\hat{\frak g}$--modules. We will show
that this category is not semisimple. Consider the following $\frak
g$--modules:
\[
M_{t}= \frac{M(t \omega _{1})}{\sum _{i=2}^{l} U(\frak
g)f_{\epsilon_{i} - \epsilon_{i+1}}v_{t \omega _{1}}},
\]
for $t \in {\Zset}_{+}$, where $M(t \omega _{1})$ denotes the Verma
$\frak g$--module with the highest weight $t \omega _{1}$ and the
highest weight vector $v_{t \omega _{1}}$. It can be easily verified
that $M_{t}$ is a highest weight $A(V_{l,1})$--module, which is not
irreducible. Using Zhu's theory we obtain a highest weight  module
for the vertex operator algebra $V_{l,1}$, which is not irreducible.
Therefore we get an example of a weak $V_{l,1}$--module from the
category $\mathcal{O}$, which is not completely reducible.
\end{remark}

\section{Vertex operator algebras associated to affine Lie algebra of
type $A_{2}^{(1)}$ with levels $n-2$, for $n \in \Nset$}
\label{sl-3}

Let $\frak g$ be the simple Lie algebra of type $A_{2}$, and
$\hat{\frak g}$ the affine Lie algebra associated to $\frak g$. In
this section we study the vertex operator algebra $N_{2}(n-2,0)$
associated to $\hat{\frak g}$ of level $n-2$, for $n \in \Nset$.

We shall extend the construction of singular vectors from \cite{A4}
in the case of the affine Lie algebra $A_2 ^{(1)}$. We will present
the formulas for singular vectors. These vectors generate ideals in
vertex operator algebras $N_2( n-2,0)$. It will be important for our
construction that top levels of these ideals are irreducible
finite-dimensional $\frak g$--modules with all $1$--dimensional
weight spaces.

\begin{theorem} \label{singular-l2}
The vector
\begin{eqnarray*}
v_{2,n}= \sum _{t=0}^{2n} \frac{1}{t!} e_{\epsilon_1 -
\epsilon_3}(-1)^t e_{\epsilon_2 - \epsilon_3}(-1)^{2n-t}
f_{\epsilon_1 - \epsilon_2}(0)^t e_{\epsilon_1 - \epsilon_2}(-1)^n
{\bf 1}
\end{eqnarray*}
is a singular vector in $N_{2}(n-2,0)$, for $n \in \Nset$.
\end{theorem}
{\bf Proof:} Since
\begin{eqnarray*}
&& e_{\epsilon_1 - \epsilon_2}(0). e_{\epsilon_1 - \epsilon_3}(-1)^t
e_{\epsilon_2 - \epsilon_3}(-1)^{2n-t} f_{\epsilon_1 -
\epsilon_2}(0)^t e_{\epsilon_1 -
\epsilon_2}(-1)^n {\bf 1} \\
&&= -(2n-t)e_{\epsilon_1 - \epsilon_3}(-1)^{t+1} e_{\epsilon_2 -
\epsilon_3}(-1)^{2n-t-1} f_{\epsilon_1 - \epsilon_2}(0)^{t}
e_{\epsilon_1 -
\epsilon_2}(-1)^n {\bf 1} \\
&&+t(2n-t+1) e_{\epsilon_1 - \epsilon_3}(-1)^t e_{\epsilon_2 -
\epsilon_3}(-1)^{2n-t} f_{\epsilon_1 - \epsilon_2}(0)^{t-1}
e_{\epsilon_1 - \epsilon_2}(-1)^n {\bf 1}
\end{eqnarray*}
for $t \in \{0,1, \ldots ,2n \}$, one can easily obtain that
$e_{\epsilon_1 - \epsilon_2}(0).v_{2,n}=0$. Similarly, one can check
that $e_{\epsilon_2 - \epsilon_3}(0).v_{2,n}=0$ and
$f_{\theta}(1).v_{2,n}=0$. \qed

Theorem \ref{singular-l2} proves that for $n \ge 2$, the maximal
submodule $N^{1}_{2}(n-2,0)$ is reducible.

Define the ideal $J_{2}(n-2,0)$ in the vertex operator algebra
$N_{2}(n-2,0)$ by
\[
J_{2}(n-2,0)=U(\hat{\frak g})v_{2,n},
\]
and denote the corresponding quotient vertex operator algebra by
$V_{2,n}$, i.e., $V_{2,n}=N_{2}(n-2,0)/J_{2}(n-2,0)$.

It follows from Zhu's theory that:

\begin{proposition} \label{prop.vekt-Zhu-sl3}
The associative algebra $A(V_{2,n})$ is isomorphic to the algebra
$U(\frak g)/I_{2,n}$, where $I_{2,n}$ is the two-sided ideal of
$U(\frak g)$ generated by
\begin{eqnarray*}
v_{2,n}'= \sum _{t=0}^{2n} \frac{1}{t!} (f_{\epsilon_1 -
\epsilon_2}^t) _L (e_{\epsilon_1 - \epsilon_2}^n) e_{\epsilon_2 -
\epsilon_3}^{2n-t} e_{\epsilon_1 - \epsilon_3}^t.
\end{eqnarray*}
\end{proposition}

The $U(\frak g)$--submodule $R$ of $U(\frak g)$ generated by the
vector $v_{2,n}'$ under the adjoint action is isomorphic to $V(3n
\omega_{2})$. Proposition \ref{one-dim} gives that  all weight
spaces of $R$ are $1$--dimensional. In particular, we have  $\dim
R_{0} =1$. The following lemma will be proved in Section
\ref{dokaz-lem}.

\begin{lemma} \label{lem.pol.sl-3}
Let
\begin{eqnarray}
&& p(h)=h_{1}(h_{1}-1)\cdot \ldots \cdot (h_{1}-n+1)
h_{2}(h_{2}-1)\cdot \ldots \cdot (h_{2}-n+1) \nonumber \\
&& \ \cdot(h_{1}+h_{2}+1)(h_{1}+h_{2})\cdot \ldots \cdot
(h_{1}+h_{2}-n+2). \label{rel.poly-l=2}
\end{eqnarray}
Then $p \in {\mathcal P}_{0}$.
\end{lemma}

Since $\dim {\mathcal P}_{0}=\dim R_{0} =1$, Corollary \ref{c.0.8}
implies that

\begin{proposition} \label{klas-zhu-A2} The set
\begin{eqnarray*}
&& \Big\{ V(t \omega _{1}+m\omega _{2}) \ \vert \ t \in \Cset,
m \in \{0,1, \ldots ,n-1 \} \Big\} \\
&& \cup \Big\{ V(m\omega _{1}+t \omega _{2}) \ \vert \ t \in \Cset,
m \in \{0,1, \ldots ,n-1 \} \Big\} \\
&& \cup \Big\{ V(t \omega _{1} + (-1-t+m) \omega _{2}) \ \vert \ t
\in \Cset, m \in \{0,1, \ldots ,n-1 \} \Big\}
\end{eqnarray*}
provides the complete list of irreducible $A(V_{2,n})$--modules from
the category $\mathcal{O}$.
\end{proposition}

In the case $n=1$, we get:

\begin{corollary}
 The set
\begin{eqnarray*}
\left\{ V(t \omega _{1}) \ \vert \ t \in {\Zset}_{+} \right\} \cup
\left\{ V(t \omega _{2}) \ \vert \ t \in {\Zset}_{+} \right\}
\end{eqnarray*}
provides the complete list of irreducible finite-dimensional
$A(V_{2,1})$--modules. Moreover,  every irreducible
finite-dimensional $A(V_{2,1})$--module has all $1$--dimensional
weight spaces for $\frak h$.
\end{corollary}

It follows from Proposition \ref{klas-zhu-A2} and Zhu's theory that:

\begin{theorem} \label{t.klas-sl3}
The set
\begin{eqnarray} && \Big\{ L(n-2,t \omega _{1}+m\omega _{2})
\ \vert \ t \in \Cset,
m \in \{0,1, \ldots ,n-1 \} \Big\} \nonumber  \\
&& \cup \Big\{ L(n-2,m\omega _{1}+t \omega _{2}) \ \vert \ t \in
\Cset,
m \in \{0,1, \ldots ,n-1 \} \Big\} \nonumber \\
&& \cup \Big\{ L(n-2,t \omega _{1} + (-1-t+m) \omega _{2}) \ \vert \
t \in \Cset, m \in \{0,1, \ldots ,n-1 \} \Big\} \label{rep-sl3}
\end{eqnarray}
provides the complete list of irreducible weak $V_{2,n}$--modules
from the category $\mathcal{O}$.
\end{theorem}

Specially, for $n=1$, we get an analogue of Theorem
\ref{t.klas.nivo-1} for the case $A_2 ^{(1)}$:

\begin{corollary}
The set
\begin{eqnarray}
 \left\{ L(-1, t \omega _{1}) \ \vert \ t
\in \Cset \right\} \cup \left\{ L(-1, t \omega _{2}) \ \vert \ t \in
\Cset \right\} \cup \left\{ L(-1, t \omega _{1}+ (-1-t) \omega _{2})
\ \vert \ t \in \Cset  \right\}  \nonumber
\end{eqnarray}
provides the complete list of irreducible weak $V_{2,1}$--modules
from the category $\mathcal{O}$.
\end{corollary}

Let $\Psi$ be the  automorphism of the vertex operator algebra
$N_2(n-2,0)$ of order two which is lifted from the automorphism
\[
e_{\alpha_1} \mapsto e_{\alpha_2}, e_{\alpha_2} \mapsto
e_{\alpha_1}, f_{\alpha_1} \mapsto f_{\alpha_2}, f_{\alpha_2}
\mapsto f_{\alpha_1}
\]
of the Lie algebra ${\frak g}= {\frak sl_3} (\Cset)$. It is also an
automorphism of the affine Lie algebra ${\hat {\frak g}}$. In a
different context, this automorphism was also considered  in
\cite{A-2007} and \cite{B-2001}.

For any ${\hat {\frak g}}$--module $M$, let $\Psi(M)$ by the ${\hat
{\frak g}}$--module obtained by applying the automorphism $\Psi$.

The proof of the following lemma is obvious.
\begin{lemma} \label{lem-o} \item[(i)]If $v$ is a non-trivial singular vector in $N_2(n-2,0)$, then $
\Psi(v)$ is also a non-trivial singular vector in $N_2(n-2,0)$.
\item[(ii)] The set (\ref{rep-sl3}) is $\Psi$--invariant.
\end{lemma}

The previous lemma implies that  the vector
\begin{eqnarray*}
\widetilde{v_{2,n}}=  \Psi(v_{2,n}) = \sum _{t=0}^{2n} \frac{(-1)
^{t}}{t!} e_{\epsilon_1 - \epsilon_3}(-1)^t e_{\epsilon_1 -
\epsilon_2}(-1)^{2n-t} f_{\epsilon_2 - \epsilon_3}(0)^t
e_{\epsilon_2 - \epsilon_3}(-1)^n {\bf 1}
\end{eqnarray*}
is a singular vector in $N_{2}(n-2,0)$.

Define the ideal $\widetilde{J_{2}}(n-2,0)= U(\hat{\frak
g})\widetilde{v_{2,n}}$ in $N_2(n-2,0)$ and denote the corresponding
quotient vertex operator algebra by $\widetilde{V_{2,n}}$.    By
using  Theorem \ref{t.klas-sl3} and Lemma \ref{lem-o} we get the
following result.

\begin{theorem}
\item[(i)] For every $n \in {\Nset}$, the vector
$\widetilde{v_{2,n}}= \Psi(v_{2,n})$ is a non-trivial singular
vector in $N_2(n-2,0)$.  In particular, the vertex operator algebra
$N_2(n-2,0)$ contains two linearly independent singular vectors of
conformal weight $3 n$.

\item[(ii)] The vertex operator algebras $V_{2,n}$ and
$\widetilde{V_{2,n}}$ are isomorphic, and $\Psi \vert _{V_{2,n}}
:V_{2,n} \rightarrow \widetilde{V_{2,n}}$  is the corresponding
isomorphism.

 \item[(iii)] The set (\ref{rep-sl3}) provides all
irreducible weak $\widetilde{V_{2,n}}$--modules from the category
${\mathcal O}$.
\end{theorem}

\begin{remark} For $n=1$ we obtain two singular vectors
in $N_2(-1,0)$  of conformal weight $3$:
\begin{eqnarray*}
&& v_{2,1}=e_{\epsilon_2 - \epsilon_3}(-1)^2 e_{\epsilon_1 -
\epsilon_2}(-1) {\bf 1} - e_{\epsilon_1 - \epsilon_3}(-1)
e_{\epsilon_2 - \epsilon_3}(-1)h_{1}(-1){\bf 1}  \\
&& \quad \ - e_{\epsilon_1 - \epsilon_3}(-1)^2 f_{\epsilon_1 -
\epsilon_2}(-1){\bf 1}, \\ && \widetilde{v_{2,1}}=e_{\epsilon_1 -
\epsilon_2}(-1)^2 e_{\epsilon_2 - \epsilon_3}(-1) {\bf 1} +
e_{\epsilon_1 - \epsilon_3}(-1)
e_{\epsilon_1 - \epsilon_2}(-1)h_{2}(-1){\bf 1}  \\
&& \quad \ - e_{\epsilon_1 - \epsilon_3}(-1)^2 f_{\epsilon_2 -
\epsilon_3}(-1){\bf 1}.
\end{eqnarray*}
Since there are no singular vectors in $N_{2}(-1,0)$ of conformal
weight $2$, we showed that the maximal ideal $N_{2}^{1}(-1,0)$ is
not generated by a single singular vector. It is an interesting
problem to describe the structure of $N_2^{1}(-1,0)$.
\end{remark}

\section{Proof of Lemma \ref{lem.pol.sl-3} } \label{dokaz-lem}

The following relations hold in $U(\frak g)$:
\begin{eqnarray}
& & (f_{\alpha}^{m}) _L (e_{\alpha}^{k}) \in (-1)^k m! {k \choose
2k-m} f_{\alpha}^{m-k} \cdot (h_{\alpha}+k-m) \cdot \ldots \cdot
(h_{\alpha} -k+1) +U(\frak g) e_{\alpha},
\nonumber \\
& &  \ \mbox{ for } 2k \geq m \geq k \mbox{ and }
 \ \forall \alpha \in \Delta_{+};  \label{rel.dokaz.1}\\
& & (f_{\alpha}^{m}) _L (e_{\alpha}^{k}) \in (-1)^m k \cdot \ldots
\cdot (k-m+1) (h_{\alpha}-k+1) \cdot \ldots \cdot
(h_{\alpha} -k+m) \cdot e_{\alpha}^{k-m}  \nonumber \\
&&  \ +U(\frak g) e_{\alpha}^{k-m+1}, \ \mbox{ for } m \leq k \mbox{
and }
 \ \forall \, \alpha \in \Delta_{+};  \label{rel.dokaz.2} \\
& & e_{\alpha}^{m}f_{\alpha}^{k} \in m \cdot \ldots \cdot (m-k+1)
(h_{\alpha}-m+k) \cdot \ldots \cdot (h_{\alpha} -m+1) \cdot
e_{\alpha}^{m-k}
+U(\frak g) e_{\alpha}^{m-k+1}, \nonumber \\
&&  \ \mbox{ for } m \geq k \mbox{ and }
 \ \forall \, \alpha \in \Delta_{+};  \label{rel.dokaz.3} \\
& & e_{\epsilon_1 - \epsilon_3}^{m}f_{\epsilon_1 - \epsilon_2 }^k
\in k \cdot \ldots \cdot (k-m+1) f_{\epsilon_1 - \epsilon_2 }^{k-m}
e_{\epsilon_2 - \epsilon_3 }^{m} +U(\frak g) e_{\epsilon_1 -
\epsilon_3 },
\ \mbox{ for } m \leq k;  \label{rel.dokaz.4} \\
& & (x^m)_L(y^k)=k \cdot \ldots \cdot (k-m+1)y^{k-m}[x,y]^m,
\ \mbox{ for } x,y \in \frak g \mbox{ such that } [[x,y],y]=0, \nonumber \\
&& \ [x,[x,y]]=0 \mbox{ and } m \leq k.  \label{rel.dokaz.5}
\end{eqnarray}

We claim that
\begin{eqnarray*}
(f_{\epsilon_1 - \epsilon_3}^n f_{\epsilon_2 - \epsilon_3}^n)_L
v_{2,n} ' \in (-1)^n (n!)^2 p(h) + U(\frak g){\frak n}_{+},
\end{eqnarray*}
with $p(h)$ given by relation (\ref{rel.poly-l=2}).

It can easily be checked that for $t>n$,
\begin{eqnarray} \label{r.2.1}
(f_{\epsilon_1 - \epsilon_3}^n f_{\epsilon_2 - \epsilon_3}^n)_L
[(f_{\epsilon_1 - \epsilon_2}^t) _L (e_{\epsilon_1 - \epsilon_2}^n)
e_{\epsilon_2 - \epsilon_3}^{2n-t} e_{\epsilon_1 - \epsilon_3}^t]
\in U(\frak g){\frak n}_{+}.
\end{eqnarray}
Let $t \leq n$. Then
\begin{eqnarray*}
&& (f_{\epsilon_1 - \epsilon_3}^n f_{\epsilon_2 - \epsilon_3}^n)_L
[(f_{\epsilon_1 - \epsilon_2}^t) _L (e_{\epsilon_1 - \epsilon_2}^n)
e_{\epsilon_2 - \epsilon_3}^{2n-t}
e_{\epsilon_1 - \epsilon_3}^t] \\
&&= {\displaystyle \sum_{{(m_{1},m_{2},m_{3}) \in {\Zset}_{+}^{3}
\atop m_{1}+m_{2}+m_{3}=n}}} {\displaystyle
\sum_{{(k_{1},k_{2},k_{3}) \in {\Zset}_{+}^{3} \atop
k_{1}+k_{2}+k_{3}=n}}} {n \choose m_{1},m_{2},m_{3}}{n \choose
k_{1},k_{2},k_{3}} (f_{\epsilon_1 - \epsilon_3}^{m_{1}}
f_{\epsilon_2 - \epsilon_3}^{k_{1}} f_{\epsilon_1 - \epsilon_2}^t)
_L (e_{\epsilon_1 -
\epsilon_2}^n) \\
&& \qquad \cdot (f_{\epsilon_1 - \epsilon_3}^{m_{2}} f_{\epsilon_2 -
\epsilon_3}^{k_{2}})_L (e_{\epsilon_1 - \epsilon_3}^t) \cdot
(f_{\epsilon_1 - \epsilon_3}^{m_{3}} f_{\epsilon_2 -
\epsilon_3}^{k_{3}})_L (e_{\epsilon_2 - \epsilon_3}^{2n-t}).
\end{eqnarray*}
Using relations (\ref{rel.dokaz.1})--(\ref{rel.dokaz.5}) one can
easily check that terms which are not in $U(\frak g){\frak n}_{+}$
are obtained only for $m_{1}+k_{1}=0$, $m_{2}+k_{2} \leq t$ and
$m_{3}+k_{3} \geq 2n-t$. Denoting $m_{2}=m$ and $k_{2}=k$, we get
\begin{eqnarray} \label{r.2.2}
&& (f_{\epsilon_1 - \epsilon_3}^n f_{\epsilon_2 - \epsilon_3}^n)_L
[(f_{\epsilon_1 - \epsilon_2}^t) _L (e_{\epsilon_1 - \epsilon_2}^n)
e_{\epsilon_2 - \epsilon_3}^{2n-t}
e_{\epsilon_1 - \epsilon_3}^t] \nonumber \\
&&\in \sum_{m=0}^{t}\sum_{k=0}^{t-m} {n \choose m} {n \choose k}
(f_{\epsilon_1 - \epsilon_2}^t) _L (e_{\epsilon_1 - \epsilon_2}^n)
\cdot (f_{\epsilon_1 - \epsilon_3}^{m} f_{\epsilon_2 -
\epsilon_3}^{k})_L
(e_{\epsilon_1 - \epsilon_3}^t)  \nonumber \\
&& \qquad \qquad \cdot (f_{\epsilon_1 - \epsilon_3}^{n-m}
f_{\epsilon_2 - \epsilon_3}^{n-k})_L (e_{\epsilon_2 -
\epsilon_3}^{2n-t}) + U(\frak g){\frak n}_{+}.
\end{eqnarray}
Using relations  (\ref{rel.dokaz.1}), (\ref{rel.dokaz.2}) and
(\ref{rel.dokaz.5}) we obtain
\begin{eqnarray*}
&& (f_{\epsilon_1 - \epsilon_3}^{n-m} f_{\epsilon_2 -
\epsilon_3}^{n-k})_L (e_{\epsilon_2 - \epsilon_3}^{2n-t})=(2n-t)
\cdot \ldots \cdot (n+m-t+1) (f_{\epsilon_2 - \epsilon_3}^{n-k})_L
(f_{\epsilon_1 - \epsilon_2}^{n-m}
e_{\epsilon_2 - \epsilon_3}^{n+m-t}) \\
&& =(2n-t) \cdot \ldots \cdot (n+m-t+1) \sum_{a=0}^{n-k} {n-k
\choose a} (f_{\epsilon_2 - \epsilon_3}^{a})_L (f_{\epsilon_1 -
\epsilon_2}^{n-m}) \cdot (f_{\epsilon_2 - \epsilon_3}^{n-k-a})_L
(e_{\epsilon_2 - \epsilon_3}^{n+m-t}) \\
&& \in (2n-t) \cdot \ldots \cdot (n+m-t+1) \sum_{a=0}^{t-k-m}
{n-k \choose a} (-1)^a (n-m)\cdot \ldots \cdot (n-m-a+1) \\
&& \cdot f_{\epsilon_1 - \epsilon_2}^{n-m-a} f_{\epsilon_1 -
\epsilon_3}^{a} \cdot (-1)^{n+m-t} (n-k-a)! {n+m-t \choose
n+2m-2t+k+a}
f_{\epsilon_2 - \epsilon_3}^{t-k-m-a} \\
&& \cdot (h_{\epsilon_2 - \epsilon_3}+m-t+k+a)\cdot \ldots \cdot
(h_{\epsilon_2 - \epsilon_3}-n-m+t+1) + U(\frak g){\frak n}_{+}.
\end{eqnarray*}
Similarly, using relations (\ref{rel.dokaz.2}) and
(\ref{rel.dokaz.5}) one can obtain that
\begin{eqnarray*}
&& (f_{\epsilon_1 - \epsilon_3}^{m} f_{\epsilon_2 -
\epsilon_3}^{k})_L
(e_{\epsilon_1 - \epsilon_3}^t) \\
&& \in t \cdot \ldots \cdot (t-k+1) \sum_{b=0}^{m}
\Big[ {m \choose b} (-1)^b b\cdot \ldots \cdot (b-k+1) \\
&& \cdot f_{\epsilon_2 - \epsilon_3}^{b} e_{\epsilon_1 -
\epsilon_2}^{k-b}
\cdot (-1)^{m-b} (t-k) \cdot \ldots \cdot (t-k-m+b+1) \\
&& \cdot (h_{\epsilon_1 - \epsilon_3}-t+k+1)\cdot \ldots \cdot
(h_{\epsilon_1 - \epsilon_3}-t+k+m-b) e_{\epsilon_1 -
\epsilon_3}^{t-k-m+b} + U(\frak g) e_{\epsilon_1 -
\epsilon_3}^{t-k-m+b+1} \Big].
\end{eqnarray*}
Using relations (\ref{rel.dokaz.3}) and (\ref{rel.dokaz.4}) one can
easily check that terms in expression \linebreak $(f_{\epsilon_1 -
\epsilon_3}^{m} f_{\epsilon_2 - \epsilon_3}^{k})_L (e_{\epsilon_1 -
\epsilon_3}^t) \cdot (f_{\epsilon_1 - \epsilon_3}^{n-m}
f_{\epsilon_2 - \epsilon_3}^{n-k})_L (e_{\epsilon_2 -
\epsilon_3}^{2n-t})$ which are not in $U(\frak g){\frak n}_{+}$ are
obtained only for $b=0$. It follows that
\begin{eqnarray} \label{r.2.3}
&& (f_{\epsilon_1 - \epsilon_3}^{m} f_{\epsilon_2 -
\epsilon_3}^{k})_L (e_{\epsilon_1 - \epsilon_3}^t) \cdot
(f_{\epsilon_1 - \epsilon_3}^{n-m} f_{\epsilon_2 -
\epsilon_3}^{n-k})_L
(e_{\epsilon_2 - \epsilon_3}^{2n-t}) \nonumber  \\
&& \in (-1)^m t\cdot \ldots \cdot (t-k-m+1) e_{\epsilon_1 -
\epsilon_2}^{k} (h_{\epsilon_1 - \epsilon_3}-t+k+1)\cdot \ldots
\cdot
(h_{\epsilon_1 - \epsilon_3}-t+k+m) \nonumber \\
&& \cdot e_{\epsilon_1 - \epsilon_3}^{t-k-m} \cdot
(2n-t) \cdot \ldots \cdot (n+m-t+1) \nonumber  \\
&& \cdot \sum_{a=0}^{t-k-m} {n-k \choose a} (-1)^a (n-m)\cdot \ldots
\cdot (n-m-a+1) f_{\epsilon_1 - \epsilon_2}^{n-m-a}
f_{\epsilon_1 - \epsilon_3}^{a} \nonumber \\
&& \cdot (-1)^{n+m-t} (n-k-a)! {n+m-t \choose n+2m-2t+k+a}
f_{\epsilon_2 - \epsilon_3}^{t-k-m-a} \nonumber  \\
&& \cdot (h_{\epsilon_2 - \epsilon_3}+m-t+k+a)\cdot \ldots \cdot
(h_{\epsilon_2 - \epsilon_3}-n-m+t+1) + U(\frak g){\frak n}_{+}.
\end{eqnarray}
It follows from relations (\ref{rel.dokaz.3}) and
(\ref{rel.dokaz.4}) that
\begin{eqnarray} \label{r.2.4}
&& e_{\epsilon_1 - \epsilon_3}^{t-k-m}f_{\epsilon_1 -
\epsilon_3}^{a} f_{\epsilon_1 - \epsilon_2}^{n-m-a}f_{\epsilon_2 -
\epsilon_3}^{t-k-m-a} \in (t-k-m) \cdot \ldots \cdot (t-k-m-a+1) \nonumber \\
&& \cdot (h_{\epsilon_1 - \epsilon_3}-t+k+m+a)\cdot \ldots \cdot
(h_{\epsilon_1 - \epsilon_3}-t+k+m+1) \nonumber \\
&& \cdot (n-m-a) \cdot \ldots \cdot (n-t+k+1) f_{\epsilon_1 -
\epsilon_2}^{n-t+k}(t-k-m+a)! \nonumber \\
&& \cdot h_{\epsilon_2 - \epsilon_3} \cdot \ldots \cdot
(h_{\epsilon_2 - \epsilon_3}-t+k+m+a+1) + U(\frak g){\frak n}_{+}.
\end{eqnarray}
Using relations  (\ref{rel.dokaz.1}), (\ref{rel.dokaz.2}),
(\ref{r.2.3}) and (\ref{r.2.4}) we get
\begin{eqnarray} \label{r.2.4.1}
&& (f_{\epsilon_1 - \epsilon_2}^t) _L (e_{\epsilon_1 -
\epsilon_2}^n) \cdot (f_{\epsilon_1 - \epsilon_3}^{m} f_{\epsilon_2
- \epsilon_3}^{k})_L (e_{\epsilon_1 - \epsilon_3}^t)
 \cdot (f_{\epsilon_1 - \epsilon_3}^{n-m} f_{\epsilon_2 - \epsilon_3}^{n-k})_L
(e_{\epsilon_2 - \epsilon_3}^{2n-t}) \nonumber \\
&& \in (-1)^n n \cdot \ldots \cdot (n-t+1) (h_{\epsilon_1 -
\epsilon_2}-n+1)
\cdot \ldots \cdot (h_{\epsilon_1 - \epsilon_2}-n+t) t! \nonumber \\
&& \cdot (h_{\epsilon_1 - \epsilon_3}-n+1)\cdot \ldots \cdot
(h_{\epsilon_1 - \epsilon_3}-n+m) \cdot (2n-t) \cdot \ldots \cdot
(n+m-t+1)(n-m)! \nonumber \\
&& \cdot h_{\epsilon_1 - \epsilon_2} \cdot \ldots \cdot
(h_{\epsilon_1 - \epsilon_2}-n+t-k+1) h_{\epsilon_2 - \epsilon_3}
\cdot \ldots \cdot (h_{\epsilon_2 - \epsilon_3}-n-m+t+1) \nonumber \\
&& \cdot \sum_{a=0}^{t-k-m} {n-k \choose a} (-1)^a
(n-k-a)! {n+m-t \choose n+2m-2t+k+a} \nonumber \\
&& \cdot (h_{\epsilon_1 - \epsilon_3}-n+m+a)\cdot \ldots \cdot
(h_{\epsilon_1 - \epsilon_3}-n+m+1)+ U(\frak g){\frak n}_{+}.
\end{eqnarray}
Clearly
\begin{eqnarray} \label{r.2.4.2}
&& \sum_{a=0}^{t-k-m} {n-k \choose a} (-1)^a
(n-k-a)! {n+m-t \choose n+2m-2t+k+a} \nonumber \\
&& \cdot (h_{\epsilon_1 - \epsilon_3}-n+m+a)\cdot \ldots \cdot
(h_{\epsilon_1 - \epsilon_3}-n+m+1) \nonumber \\
&& =(n-k)! \sum_{a=0}^{t-k-m} {-h_{\epsilon_1 - \epsilon_3}+n-m-1
\choose a}
{n+m-t \choose t-k-m-a} \nonumber \\
&& =(n-k)! {-h_{\epsilon_1 - \epsilon_3}+2n-t-1 \choose t-k-m}.
\end{eqnarray}
It follows from relations (\ref{r.2.2}), (\ref{r.2.4.1}) and
(\ref{r.2.4.2}) that
\begin{eqnarray*}
&& (f_{\epsilon_1 - \epsilon_3}^n f_{\epsilon_2 - \epsilon_3}^n)_L
[(f_{\epsilon_1 - \epsilon_2}^t) _L (e_{\epsilon_1 - \epsilon_2}^n)
e_{\epsilon_2 - \epsilon_3}^{2n-t}
e_{\epsilon_1 - \epsilon_3}^t] \\
&& \in (-1)^n n \cdot \ldots \cdot (n-t+1)t! (h_{\epsilon_1 -
\epsilon_2}-n+1)
\cdot \ldots \cdot (h_{\epsilon_1 - \epsilon_2}-n+t) \\
&& \cdot \sum_{m=0}^{t} {n \choose m} (h_{\epsilon_1 -
\epsilon_3}-n+1)\cdot \ldots \cdot
(h_{\epsilon_1 - \epsilon_3}-n+m) \\
&& \cdot (2n-t) \cdot \ldots \cdot (n+m-t+1)(n-m)! h_{\epsilon_2 -
\epsilon_3}
\cdot \ldots \cdot (h_{\epsilon_2 - \epsilon_3}-n-m+t+1) \\
&& \cdot \sum_{k=0}^{t-m} {n \choose k} h_{\epsilon_1 - \epsilon_2}
\cdot \ldots \cdot (h_{\epsilon_1 - \epsilon_2}-n+t-k+1) \\
&& \cdot (n-k)! {-h_{\epsilon_1 - \epsilon_3}+2n-t-1 \choose t-k-m}+
U(\frak g){\frak n}_{+}.
\end{eqnarray*}
Furthermore
\begin{eqnarray*}
&&  \sum_{k=0}^{t-m} {n \choose k} h_{\epsilon_1 - \epsilon_2} \cdot
\ldots \cdot (h_{\epsilon_1 - \epsilon_2}-n+t-k+1) (n-k)!
{-h_{\epsilon_1 - \epsilon_3}+2n-t-1 \choose
t-k-m} \\
&& = n! h_{\epsilon_1 - \epsilon_2} \cdot \ldots \cdot
(h_{\epsilon_1 - \epsilon_2}-n+t+1) \sum_{k=0}^{t-m} {h_{\epsilon_1
- \epsilon_2}-n+t \choose k} {-h_{\epsilon_1 - \epsilon_3}+2n-t-1
\choose
t-m-k} \\
&& = n! h_{\epsilon_1 - \epsilon_2} \cdot \ldots \cdot
(h_{\epsilon_1 - \epsilon_2}-n+t+1) {-h_{\epsilon_2 -
\epsilon_3}+n-1 \choose t-m},
\end{eqnarray*}
which implies
\begin{eqnarray} \label{r.2.5}
&& (f_{\epsilon_1 - \epsilon_3}^n f_{\epsilon_2 - \epsilon_3}^n)_L
[(f_{\epsilon_1 - \epsilon_2}^t) _L (e_{\epsilon_1 - \epsilon_2}^n)
e_{\epsilon_2 - \epsilon_3}^{2n-t}
e_{\epsilon_1 - \epsilon_3}^t] \nonumber \\
&& \in (-1)^n n! \cdot n \cdot \ldots \cdot (n-t+1)t! \cdot
h_{\epsilon_1 - \epsilon_2}
\cdot \ldots \cdot (h_{\epsilon_1 - \epsilon_2}-n+1) \nonumber \\
&& \cdot h_{\epsilon_2 - \epsilon_3} \cdot \ldots \cdot
(h_{\epsilon_2 - \epsilon_3}-n+1)
\sum_{m=0}^{t} {n \choose m}(n-m)! (-1)^{t-m} \frac{1}{(t-m)!} \nonumber \\
&& \cdot (h_{\epsilon_1 - \epsilon_3}-n+1)\cdot \ldots \cdot
(h_{\epsilon_1 - \epsilon_3}-n+m) (2n-t) \cdot \ldots \cdot
(n+m-t+1)+ U(\frak g){\frak n}_{+}. \nonumber \\
&&
\end{eqnarray}
One can easily verify that
\begin{eqnarray} \label{r.2.6}
&& \sum_{m=0}^{t} {n \choose m}(n-m)! (-1)^{t-m} \frac{1}{(t-m)!}
(h_{\epsilon_1 - \epsilon_3}-n+1)\cdot \ldots \cdot
(h_{\epsilon_1 - \epsilon_3}-n+m) \nonumber \\
&& \cdot (2n-t) \cdot \ldots \cdot (n+m-t+1) \nonumber \\
&& = n! \cdot (2n-t) \cdot \ldots \cdot (n+1) (-1)^t {-h_{\epsilon_1
- \epsilon_3}+2n-1 \choose t}.
\end{eqnarray}
It follows from Proposition \ref{prop.vekt-Zhu-sl3} and relations
(\ref{r.2.1}), (\ref{r.2.5}) and (\ref{r.2.6}) that
\begin{eqnarray*}
&& (f_{\epsilon_1 - \epsilon_3}^n f_{\epsilon_2 - \epsilon_3}^n)_L
v_{2,n}' \in (-1)^n (n!)^2 \cdot h_{\epsilon_1 - \epsilon_2}
\cdot \ldots \cdot (h_{\epsilon_1 - \epsilon_2}-n+1) \\
&& \cdot h_{\epsilon_2 - \epsilon_3} \cdot \ldots \cdot
(h_{\epsilon_2 - \epsilon_3}-n+1)
\sum_{t=0}^{n} (2n-t) \cdot \ldots \cdot (n-t+1) \\
&& \cdot (-1)^t {-h_{\epsilon_1 - \epsilon_3}+2n-1 \choose t}.
\end{eqnarray*}
Finally
\begin{eqnarray*}
&& \sum_{t=0}^{n} (2n-t) \cdot \ldots \cdot (n-t+1)
(-1)^t {-h_{\epsilon_1 - \epsilon_3}+2n-1 \choose t} \\
&& = n! \cdot (-1)^n {-h_{\epsilon_1 - \epsilon_3}+n-2 \choose n},
\end{eqnarray*}
which implies
\begin{eqnarray*}
&& (f_{\epsilon_1 - \epsilon_3}^n f_{\epsilon_2 - \epsilon_3}^n)_L
v_{2,n}' \in (-1)^n (n!)^2 \cdot h_{\epsilon_1 - \epsilon_2} \cdot
\ldots \cdot (h_{\epsilon_1 - \epsilon_2}-n+1) \\ && \cdot
h_{\epsilon_2 - \epsilon_3} \cdot \ldots \cdot (h_{\epsilon_2 -
\epsilon_3}-n+1) (h_{\epsilon_1 - \epsilon_3}+1) \cdot \ldots \cdot
(h_{\epsilon_1 - \epsilon_3}-n+2),
\end{eqnarray*}
and the proof is complete.

\section{Conclusions}

In this paper we constructed two families of isomorphic operator
algebras $V_{2,n}$ and $\widetilde{V_{2,n}}$ associated to affine
Lie algebra of type $A_2 ^{(1)}$ with positive integer levels $n-2$,
for $n \ge 2$. These vertex operator algebras are non-simple and
they are different both from $N_{2}(n-2,0)$ and its simple quotient
$L_{2}(n-2,0)$ which were previously extensively studied in
\cite{DLM,FZ,KWn,L,MP}. The class of irreducible weak
$N_{2}(n-2,0)$--modules includes all irreducible $A_2
^{(1)}$--modules of level $n-2$ from the category $\mathcal{O}$.
Vertex operator algebra $L_{2}(n-2,0)$ has finitely many irreducible
weak modules and any irreducible weak $L_{2}(n-2,0)$--module is
integrable highest weight $A_2 ^{(1)}$--module of level $n-2$.
Moreover, the category of weak $L_{2}(n-2,0)$--modules is
semisimple. By studying the representation theory of $V_{2,n}$, we
obtained a new interesting subcategory of the category of $A_2
^{(1)}$--modules of level $n-2$ from the category $\mathcal{O}$. The
irreducible objects in that category include irreducible weak
$L_{2}(n-2,0)$--modules, but this class is smaller than the class of
irreducible weak $N_{2}(n-2,0)$--modules from the category
$\mathcal{O}$. Irreducible weak $V_{2,n}$--modules from the category
$\mathcal{O}$ are parameterized by certain lines in $\hat {\frak
h}^{*}$, which implies that there are infinitely many of them.
Furthermore, the category of  $V_{2,n}$--modules  is not semisimple,
which makes the representation theory of $V_{2,n}$ much different
from the representation theory of $L_{2}(n-2,0)$.

We also gave examples of vertex operator algebras associated to
affine Lie algebras whose level is neither a positive integer nor an
admissible rational number. We studied representation theory of
vertex operator algebras $V_{l,1}$ associated to affine Lie algebras
of type $A_l ^{(1)}$ with negative integer level $-1$, for $l \ge
2$. These vertex operator algebras are quotients of generalized
Verma modules $N_{l}(-1,0)$ by certain ideals defined in \cite{A4}
and in this paper. In the case $l=2$, we constructed two isomorphic
vertex operator algebras $V_{2,1}$ and $\widetilde{V_{2,1}}$. We
showed that irreducible weak $V_{l,1}$--modules from category
$\mathcal{O}$ are parameterized by certain lines in $\hat {\frak
h}^{*}$, and that the category of weak $V_{l,1}$--modules  is not
semisimple.

In  our future work we plan to study fusion rules for the
irreducible representations classified in this paper.

\end{document}